\documentclass[a4paper,12pt]{article}
\usepackage{amsmath,amssymb}
\usepackage{color}
\usepackage{epsfig}
\usepackage{graphicx}

\definecolor{green1}{rgb}{0.1,0.65,0}
\definecolor{green2}{rgb}{0.1,0.6,0.1}
\definecolor{blue1}{rgb}{0.14,0.6,1.0}
\definecolor{blue3}{rgb}{0.05,0.05,0.5}
\definecolor{viol}{rgb}{0.4,0,0.9}

\definecolor{cola}{rgb}{0.4,0,0.4}
\definecolor{colb}{rgb}{0,0.3,0.4}
\definecolor{colc}{rgb}{0.4,0.3,0}
\definecolor{oran}{rgb}{.95,.35,0}
\definecolor{brass}{rgb}{0.6,0.15,0}
\definecolor{whit}{rgb}{1,1,1}
\definecolor{grey1}{rgb}{.6,.6,.6}
\definecolor{black}{rgb}{0,0,0}

\newtheorem{theorem}{Theorem}  [section]
\newtheorem{proposition}[theorem]{Proposition}
\newtheorem{lemma}[theorem]{Lemma}

\newtheorem{remark}[theorem]{Remark}

\numberwithin{equation}{section}

\font \lgr cmmib10 scaled \magstep1

\newfont{\bcal}{cmbsy10 scaled \magstep1}
\newfont{\ccal}{cmsy10 scaled \magstep1}
\newfont{\ctv}{msam10}

\newcommand{\bbox}{\mbox{\ctv \symbol{4}}}
\def\QED{{$\hfill\bbox$}}
\newenvironment{pf}[1]{\par\vskip1mm{\noindent\it #1.}\ }{\QED\par\vskip2mm}

\newcommand{\real}{\mathbb{R}}

\newcommand{\tens}{\mathbb{T}_{\rm sym}^{3\times 3}}

\def\dive{\mbox{\rm div\,}}
\def\sign{\mbox{\rm sign}}

\def\supess{\mathop{\mbox{sup\,ess}\,}}

\def\expe{{\rm e}}

\def\mae{\mbox{a.e.}}

\def\vect#1#2{\left(\!\begin{array}{ll} #1\\#2\end{array}\right)\!}

\def\bfsi{\mbox{\lgr \char27}}
\def\bfe{\mbox{\lgr \char34}}

\def\bfde{\mbox{\lgr \char14}}

\def\tbfe{\tilde\bfe}

\newcommand{\io}{\int_\Omega}
\newcommand{\ipo}{\int_{\partial\Omega}}

\def\bfA{\mathbf{A}}
\def\bfb{\mathbf{b}}

\def\bfg{\mathbf{g}}
\def\bfk{\mathbf{k}}
\def\bfn{\mathbf{n}}
\def\bfq{\mathbf{q}}
\def\bfu{\mathbf{u}}
\def\bfv{\mathbf{v}}

\def\bfdd{\mbox{\bf \,:\,}}

\def\dd{\mathrm{\,d}}

\def\barr{\begin{array}}
\def\earr{\end{array}}

\def\bpf{\begin{pf}}
\def\epf{\end{pf}}

\numberwithin{equation}{section}

\begin{document}

\title{\bf A bottle in a freezer
\footnote{Supported by the DFG Research Center {\sc Matheon}}}

\author{
Pavel Krej\v{c}\'{\i}
\footnote{Weierstrass Institute for Applied
Analysis and Stochastics, Mohrenstr.~39, D-10117 Berlin,
Germany, and Mathematical Institute,
Czech Academy of Sciences, \v{Z}itn\'a 25, CZ-11567 Praha 1, Czech Republic,
E-mail {\tt  krejci@wias-berlin.de, krejci@math.cas.cz}},
Elisabetta Rocca
\footnote{Dipartimento di Matematica, Universit\`a di Milano,
Via Saldini 50, 20133 Milano, Italy, E-mail {\tt elisabetta.rocca@unimi.it}}
\footnote{A large part of this work was done during E.~Rocca's visit
at WIAS Berlin in October/November 2008.},
and J\"urgen Sprekels
\footnote{Weierstrass Institute for Applied
Analysis and Stochastics, Mohrenstr.~39, D-10117 Berlin,
Germany, E-mail {\tt  sprekels@wias-berlin.de}}}
\date{}
\maketitle
\date{}


\vspace{-.4cm}

\noindent {\bf Abstract.}
We propose here a model for solidification of a liquid contents
of an elastic bottle in a freezer. The main goal is to explain
the occurrence of high stresses inside the bottle. As a by-product,
we derive a formula for the undercooling coefficient
in terms of the elasticity constants, latent heat, and the
phase expansion coefficient.
We investigate the well-posedness of the three-dimensional model:
we prove the existence and uniqueness of a solution for the corresponding
initial-boundary value problem which couples a PDE with an integrodifferential equation and
an ordinary differential inclusion ruling the evolution of the phase parameter.
Finally, we prove some results on the long time behavior of solutions.

\vspace{.4cm}

\vspace{4mm}


\setcounter{table}{0}

\vspace{1mm}

\section{Introduction}\label{intr}

We derive a simple model for solid-liquid phase transition of a medium inside
an elastic container. The main goal is to give a qualitative and quantitative
description of the interaction between volume, pressure, phase, and temperature
changes in the situation that the specific volume of the solid phase exceeds
the specific volume of the liquid phase. We compute
the undercooling coefficient for the special case of water and ice.

There is an abundant classical literature on the study of phase transition processes, see e.g.
the monographs \cite{bs}, \cite{fremond}, \cite{visintin} and the references therein.
In \cite{fr1}, the authors proposed to interpret a phase transition process in terms of a
balance equation for macroscopic motions, and to include the possibility of voids.
Well-posedness of an initial-boundary value problem associated with the resulting
PDE system is proved there.

The microscopic approach has been pursued in \cite{fr2} in the case of two different
densities $\varrho_1$ and $\varrho_2$ for the two substances undergoing phase transitions.
The evolution of a liquid substance, e.g., water, in a rigid container subject to freezing
is described by a mass balance in Eulerian coordinates, an entropy balance, and a phase field
equation. The flow is governed by a counterpart of the Darcy law.
Since the density $\varrho_{2}$ of ice is lower
than the density $\varrho_{1}$ of water, experiments -- for instance the freezing of a glass
bottle filled with water -- show that the water pressure increases
up to the rupture of the bottle. When the container is not
impermeable, freezing may produce a non-homogeneous material, for
instance water ice or sorbet. This particular example is treated in \cite{fr2} where the model
is presented and a suitable variational formulation  of the resulting nonlinear and singular
PDE system is solved. In the present paper, we have also other applications in mind.

Let us also mention the papers \cite{rr1} and \cite{rr2} dealing
with macroscopic stresses in phase transitions models, where the
different properties of the viscous (liquid) and elastic (solid)
phases are taken into account and the coexisting viscous and
elastic properties of the system are given a distinguished role,
under the working assumption that they indeed influence the phase
transition process. The model there includes inertia,
viscous, and shear viscosity effects (depending on the phases),
while thermal and phase expansion of the substance are
neglected. This is reflected in the analytical
expressions of the associated PDEs for the strain $\bfu$ and
the phase parameter $\chi$: the $\chi$-dependence, e.g.,  in the
stress-strain relation leads to the possible degeneracy of the
elliptic operator therein.  In \cite{rr1} and \cite{rr2},
respectively, local existence (in the 3D case) and well-posedness
(in the 1D case) for the corresponding initial-boundary value
problems are proved. Finally, we can quote in this framework the
model analyzed  in~\cite{kss2} and \cite{kss}, which
pertains to nonlinear thermoviscoplasticity: in the
one-dimensional (in space) case, the authors prove the global
well-posedness of a PDE system, incorporating both hysteresis
effects and modeling phase change, which however does not display
a degenerating character.

Here, in Section~\ref{mode}, we derive a completely different model without referring
to any microscopic balance laws, and deal exclusively with physically measurable quantities.
We assume that the displacements are small. This enables us to state the system in
Lagrangian coordinates. The main difference with respect to the Eulerian framework in \cite{fr2}
is that in Lagrangian coordinates,
the mass conservation law is equivalent to the same constant mass density in liquid and in solid,
but the specific volumes of the liquid and solid phases are different.
For simplicity, we assume that the speed of sound, specific heat, heat conductivity,
viscosity, and thermal expansion coefficient do not depend on the phase,
the evolution is slow, and the shear viscosity, shear stresses, and inertia effects
are negligible. The process is driven by energy balance, quasistatic
momentum balance, and a phase dynamics equation.
Still in Section~\ref{mode}, we verify the thermodynamic
consistency of the model, and in Section \ref{equi} we study the equilibria.
We observe there that a pure solid state can only be reached
if the external temperature is below a certain threshold, which is lower
than the freezing point and depends in particular on the elasticity of the boundary.
For water and ice, we explicitly compute the undercooling rate,
which turns out to be around $5\%$ if the container is rigid.
For intermediate temperatures between freezing point and undercooling limit,
there exists a continuum of distinct equilibria with mixtures of solid and liquid.
If in this situation the bottle breaks, an instantaneous solidification
takes place.

The well-posedness of the three-dimensional model is investigated in Section \ref{exiu},
and the asymptotic stabilization of the process is proved in Section \ref{long}.

\section{The model}\label{mode}

As reference state, we consider a liquid substance contained in a bounded connected
bottle $\Omega \subset \real^3$ with boundary of class $C^{1,1}$. The state variables are
the absolute temperature $\theta>0$, the displacement $\bfu \in \real^3$,
and the phase variable $\chi \in [0,1]$. The value  $\chi = 0$
means solid, $\chi = 1$ means liquid, $\chi \in (0,1)$ is a mixture
of the two.

We make the following modeling hypotheses.

\begin{itemize}
\item[{\bf (A1)}]
The displacements are small.
Therefore, we state the problem in {\em Lagrangian coordinates\/},
in which the mass conservation is equivalent to the condition
of a constant mass density $\varrho_0>0$.
\item[{\bf (A2)}]
The substance is compressible, and the speed of sound does not depend on the phase.
\item[{\bf (A3)}]
The evolution is slow, and we neglect shear viscosity and inertia effects.
\item[{\bf (A4)}]
We neglect shear stresses and gravity effects.
\end{itemize}

In agreement with {\bf (A1)}, we define the strain $\bfe$ as an element
of the space $\tens$ of symmetric tensors by the formula
\begin{equation}\label{eps}
\bfe = \nabla_s \bfu := \frac12 (\nabla \bfu + (\nabla \bfu)^T).
\end{equation}
Let $\bfde\in \tens$ denote the Kronecker tensor.
By {\bf (A4)}, the elasticity matrix $\bfA$ has the form
\begin{equation}\label{ela}
\bfA\bfe = \lambda (\bfe\bfdd\bfde)\,\bfde\,,
\end{equation}
where ``$\bfdd$'' is the canonical scalar product in $\tens$,
and $\lambda > 0$ is the Lam\'e constant (or {\em bulk elasticity modulus}),
which we assume to be independent of $\chi$ by virtue of {\bf (A2)}.
Note that $\lambda$ is related to the speed of sound $v_0$ by the formula
$v_0 = \sqrt{\lambda/\varrho_0}$.

We want to model the situation where the specific volume $V_{solid}$
of the solid phase is larger than the specific volume $V_{liquid}$ of the
liquid phase. Considering the liquid phase as the reference state, we
{introduce} the dimensionless phase expansion coefficient
$\alpha = (V_{solid}-V_{liquid})/V_{liquid} > 0$,
and we define the phase expansion strain $\tbfe$ by
\begin{equation}\label{tbfe}
\tbfe(\chi) = \frac{\alpha}{3}(1-\chi) \bfde\,.
\end{equation}
We fix positive constants $c_0$ (specific heat), $L_0$ (latent heat),
$\theta_c$ (freezing point at standard atmospheric pressure), $\gamma_0$
(phase relaxation coefficient), $\beta$ (thermal expansion coefficient), and consider the
specific free energy $f$ in the form
\begin{eqnarray}\label{free}
f &=& c_0\theta\Big(1-\log\Big(\frac{\theta}{\theta_c}\Big)\Big) + \frac{\lambda}{2\varrho_0}
((\bfe - \tbfe(\chi))\bfdd\bfde)^2 - \frac{\beta}{\varrho_0}(\theta - \theta_c)\bfe\bfdd\bfde
\\[2mm]\nonumber
&& +\, L_0\left(\chi\left(1 - \frac{\theta}{\theta_c}\right)
+  I(\chi)\right)\,,
\end{eqnarray}
where $I$ is the indicator function of the interval $[0,1]$.

To derive the balance equations, we first proceed formally,
assuming that the temperature is positive.
This assumption will be justified in the subsequent sections.
The stress tensor $\bfsi$ is decomposed into the sum $\bfsi^v + \bfsi^e$
of the viscous component  $\bfsi^v$ and elastic component $\bfsi^e$.
The state functions $\bfsi^v$,$\bfsi^e$, $s$ (specific entropy),
and $e$ (specific internal energy) are given by the formulas
\begin{eqnarray}\label{e1}
\bfsi^v &=& \nu (\bfe_t\bfdd\bfde)\bfde\\[2mm]\label{e2a}
\bfsi^e &=& \varrho_0\frac{\partial f}{\partial \bfe} =
\left(\lambda (\bfe\bfdd\bfde - \alpha(1-\chi)) - \beta (\theta - \theta_c) \right)
\bfde \,,\\[2mm]\label{e2}
s &=& -\frac{\partial f}{\partial \theta} =
c_0\log\left(\frac{\theta}{\theta_c}\right)  + \frac{L_0}{\theta_c}\chi
+ \frac{\beta}{\varrho_0} \bfe\bfdd\bfde \,,\\[2mm]\label{e3}
e &=& f + \theta\,s = c_0\theta  + \frac{\lambda}{2\varrho_0} (\bfe\bfdd\bfde - \alpha(1-\chi))^2
+ \frac{\beta}{\varrho_0}\theta_c \bfe\bfdd\bfde + L_0(\chi + I(\chi))\,,
\end{eqnarray}
where $\nu>0$ is the volume viscosity coefficient. The scalar quantity
\begin{equation}\label{press}
p := -\nu \bfe_t\bfdd\bfde - \lambda (\bfe\bfdd\bfde - \alpha(1-\chi))
+ \beta (\theta - \theta_c)
\end{equation}
is the {\em pressure\/} and the stress has the form $\bfsi = -p\,\bfde$.
The process is governed by the balance equations
\begin{eqnarray}
\dive \bfsi &=& 0 \hspace{17mm}  \mbox{(mechanical equilibrium)}\label{bal1}\\[2mm]
\varrho_0 e_t + \dive \bfq &=& \bfsi\bfdd\bfe_t \hspace{10mm}  \mbox{(energy balance)}\label{bal2} \\[2mm]
-\gamma_0\chi_t &\in& \partial_\chi f \hspace{12mm}  \mbox{(phase relaxation law)}\label{bal3}
\end{eqnarray}
where $\partial_\chi$ is the partial subdifferential with respect to $\chi$,
and $\bfq$ is the heat flux vector that we assume in the form
\begin{equation}\label{flux}
\bfq = -\kappa\nabla\theta
\end{equation}
with a constant heat conductivity $\kappa > 0$.
The equilibrium equation (\ref{bal1}) can be rewritten in the form $\nabla p = 0$, hence
\begin{equation}\label{bala1}
p(x,t) = P_{stand} + P(t)\,,
\end{equation}
where $P_{stand}$ is the constant standard pressure, and $P$ is
a function of time only, which is to be determined. We assume the external pressure
in the form $P_{ext} = P_{stand} + p_0$ with a constant deviation $p_0$.
The normal force acting on the boundary is $-(\bfsi + P_{ext}\bfde) \bfn
= (P(t) -p_0)\bfde \bfn = (P(t) -p_0) \bfn$, where $\bfn$ denotes the unit outward normal vector
(this notation is slightly ambiguous: the first two terms in this vector
identity involve left multiplication of a vector by a matrix, while the last term is a vector
multiplied by a scalar).
We assume an elastic response of the boundary, and a heat transfer
proportional to the inner and outer temperature difference. On $\partial \Omega$,
we thus prescribe boundary conditions for $\bfu$ and $\theta$ in the form
\begin{eqnarray}\label{bcu}
(P(t) -p_0)\bfn &=& \bfk(x)\bfu \,,\\[2mm]\label{bctheta}
\bfq\cdot \bfn &=& h(x) (\theta - \theta_\Gamma)
\end{eqnarray}
with a given symmetric positive definite matrix $\bfk$ (elasticity of the boundary),
a positive function $h$ (heat transfer coefficient), and a constant $\theta_\Gamma >0$ (external
temperature).
This enables us to find an explicit relation between $\dive \bfu$ and $P$.
Indeed, on $\partial \Omega$ we have by (\ref{bcu}) that
$\bfu\cdot\bfn = (P(t) - p_0)\bfk^{-1}(x)\bfn(x)\cdot\bfn(x)$.
Assuming that $\bfk^{-1}\bfn\cdot\bfn$ belongs to $L^1(\partial \Omega)$, we set
\begin{equation}\label{bck}
\frac{1}{K_{\Gamma}} \ = \ \ipo \bfk^{-1}(x)\bfn(x)\cdot\bfn(x) \dd s(x)\,,
\end{equation}
and obtain by Gauss' Theorem that
\begin{equation}\label{bcu2}
U_\Omega(t) := \io \dive \bfu(x,t) \dd x = \frac{1}{K_{\Gamma}}(P(t) - p_0)\,.
\end{equation}
Under the small strain hypothesis, the function $\dive \bfu$ describes the
local relative volume increment. Hence, Eq.~(\ref{bcu2}) establishes a linear
relation between the total relative volume increment $U_\Omega(t)$ and the relative pressure
$P(t) - p_0$.
We have $\bfe\bfdd\bfde = \dive \bfu$, and thus the mechanical equilibrium equation
(\ref{bala1}), due to \eqref{press} and \eqref{bcu2}, reads
\begin{equation}\label{equau}
\nu \dive \bfu_t + \lambda (\dive \bfu - \alpha(1-\chi)) - \beta
(\theta - \theta_c) = -p_0 - K_{\Gamma}U_\Omega(t)\,.
\end{equation}
As a consequence of (\ref{free}), the energy balance and the phase relaxation equation
in (\ref{bal2})--(\ref{bal3}) have the form
\begin{eqnarray}\label{equ1}
\varrho_0 c_0\theta_t - \kappa\Delta\theta &=& \nu (\dive \bfu_t)^2 - \beta\theta \dive \bfu_t -
\big(\alpha\lambda(\dive \bfu - \alpha(1-\chi)) + \varrho_0 L_0\big)\chi_t\,,\qquad\\ [2mm]\label{equ4}
 -\varrho_0\gamma_0 \chi_t &\in& \alpha\lambda
(\dive\bfu - \alpha(1-\chi)) + \varrho_0 L_0 \left(1 - \frac{\theta}{\theta_c}
+ \partial I(\chi)\right)\,,
\end{eqnarray}
where $\partial$ denotes the subdifferential. For simplicity, we now set
\begin{equation}\label{zero}
c := \varrho_0 c_0\,, \quad \gamma := \varrho_0\gamma_0\,, \quad L := \varrho_0 L_0\,.
\end{equation}
The system now completely decouples. For the unknown functions $\theta, \chi$, and $U = \dive\bfu$,
we have a closed system of one PDE and two ``ODEs''
(note that mathematically, $\partial I(\chi)$ is the same as $L \partial I(\chi)$)
\begin{eqnarray}\label{sys2}
c\theta_t - \kappa\Delta\theta &=& \nu U_t^2 - \beta\theta U_t -
\big(\alpha\lambda(U - \alpha(1-\chi)) + L\big)\chi_t\,,\\ [2mm]\label{sys1}
\nu U_t + \lambda U &=&
\alpha\lambda (1-\chi) + \beta (\theta - \theta_c)
- p_0 - K_{\Gamma}U_\Omega(t)\,,\\ [2mm]\label{sys3}
-\gamma \chi_t &\in& \alpha\lambda (U - \alpha(1-\chi))
+ L \left(1 - \frac{\theta}{\theta_c}\right) + \partial I(\chi)\,,
\end{eqnarray}
with $U_\Omega(t) = \ipo U(x,t) \dd s(x)$, and
with boundary condition (\ref{bctheta}), (\ref{flux}).
To find $\bfu$, we first define $\Phi$ as a solution of the Poisson
equation $\Delta\Phi = U$ with the Neumann boundary condition
$\nabla\Phi\cdot \bfn = K_{\Gamma} U_\Omega(t) \bfk^{-1}(x)\bfn(x)\cdot\bfn(x)$.
With this $\Phi$, we find $\tilde \bfu$ as a solution to the problem
\begin{eqnarray}\label{sys4a}
\dive \tilde\bfu &=& 0\hspace{7mm} \mbox{in }\ \Omega\times (0,\infty)\,,\\ \label{sys6a}
\tilde\bfu \cdot \bfn \ = \ 0\,, \quad
(\tilde\bfu + \nabla\Phi - K_{\Gamma} U_\Omega \bfk^{-1}\bfn)\times\bfn &=& 0\hspace{7mm}
\mbox{on }\ \partial\Omega\times (0,\infty)\,,
\end{eqnarray}
and set $\bfu = \tilde\bfu + \nabla\Phi$. Then $\bfu$ satisfies a.e. in $\Omega$ the equation
$\dive \bfu = U$, together with the boundary condition
(\ref{bcu}), that is, $\bfu =K_{\Gamma} U_\Omega \bfk^{-1}\bfn$ on $\partial\Omega$.

For the solution to (\ref{sys4a})--(\ref{sys6a}), we refer to
\cite[Lemma 2.2]{gr} which states that
for each $\bfg\in H^{1/2}(\partial\Omega)^3$ satisfying $\ipo \bfg\cdot\bfn \dd s(x) = 0$
there exists a function $\tilde\bfu\in H^1(\Omega)^3$, unique up to an additive
function $\bfv$ from the set $V$ of divergence-free $H^1(\Omega)$ functions vanishing on
$\partial\Omega$, such that $\dive \tilde\bfu = 0$ in $\Omega$, $\tilde\bfu = \bfg$
on $\partial\Omega$. In terms of the system (\ref{sys4a})--(\ref{sys6a}), it suffices to set
$\bfg =  ((\nabla\Phi - K_{\Gamma} U_\Omega \bfk^{-1}\bfn)\times\bfn)\times\bfn$
and use the identity $(\bfb\times\bfn)\times\bfn = (\bfb\cdot\bfn)\,\bfn-\bfb$ for every
vector $\bfb$. Moreover, the estimate
\begin{equation}\label{turbu}
\inf_{\bfv \in V}\|\tilde \bfu + \bfv\|_{H^1(\Omega)} \le C\,\|\bfg\|_{H^{1/2}(\partial\Omega)}
\le \tilde C \|\Phi\|_{H^2(\Omega)}
\end{equation}
holds with some constants $C, \tilde C$. The required regularity is available here
by virtue of the assumption that $\Omega$ is of class $C^{1,1}$, provided $\bfk^{-1}$ belongs
to $H^{1/2}(\partial \Omega)$. Note that a weaker formulation
of problem (\ref{sys4a})--(\ref{sys6a}) can be found in \cite[Section 4]{ag}.

Due to our hypotheses {\bf (A3)}, {\bf (A4)}, we thus lose any control on
possible volume preserving turbulences $\bfv \in V$. This, however, has no influence
on the system (\ref{sys2})--(\ref{sys3}), which is the subject of our interest here.
Inequality (\ref{turbu}) shows that $\bfv \in V$ can be chosen in such a way that
hypothesis {\bf (A1)} is not violated.

In terms of the new variables $\theta, U, \chi$, the energy $e$ and entropy $s$ can be written~as
\begin{eqnarray}\label{energy}
e &=& c_0\theta + \frac{\lambda}{2\varrho_0}(U - \alpha(1-\chi))^2
+ \frac{\beta}{\varrho_0} \theta_c U
+ L_0(\chi + I(\chi))\,,\\[2mm]\label{entropy}
s &=& c_0\log\left(\frac{\theta}{\theta_c}\right)  + \frac{L_0}{\theta_c}\chi
+ \frac{\beta}{\varrho_0} U \,.
\end{eqnarray}
The energy functional has to be supplemented with the boundary
energy term
\begin{equation}\label{enerb}
E_\Gamma(t) \ = \ \frac{K_{\Gamma}}{2}\left(U_{\Omega}(t) +\frac{p_0}{K_{\Gamma}}\right)^2\,.
\end{equation}
The energy and entropy balance equations now read
\begin{eqnarray}\label{princ1}
\frac{\dd}{\dd t} \left(\io \varrho_0 e(x,t)\dd x
+ E_\Gamma(t) \right)
&=& \ipo h(x)(\theta_\Gamma - \theta)\dd s(x)\,,\\[2mm]\label{princ2}
\varrho_0 s_t + \dive \frac{\bfq}{\theta} &=& \frac{\kappa |\nabla\theta|^2}{\theta^2}
+ \frac{\gamma}{\theta} \chi_t^2 + \frac{\nu}{\theta}U_t^2 \ \ge\ 0\,,\\[2mm]\label{princ3}
\frac{\dd}{\dd t}\io \varrho_0 s(x,t)\dd x &=& \ipo \frac{h(x)}{\theta}(\theta_\Gamma - \theta)
\dd s(x) \\[1mm]\nonumber
&& +\, \io\left(\frac{\kappa |\nabla\theta|^2}{\theta^2}
+ \frac{\gamma}{\theta} \chi_t^2 + \frac{\nu}{\theta}U_t^2  \right)\dd x\,.
\end{eqnarray}
The entropy balance (\ref{princ2}) says that the entropy production
on the right hand side is nonnegative in agreement with the second principle
of thermodynamics. The system is not closed, and the energy supply through
the boundary is given by the right hand side of (\ref{princ1}).

We prescribe the initial conditions
\begin{eqnarray}\label{ini2}
\theta(x,0) &=& \theta^0(x)\\ \label{ini1}
U(x,0) &=& U^0(x)\\ \label{ini3}
\chi(x,0) &=& \chi^0(x)
\end{eqnarray}
for $x \in \Omega$, and compute from (\ref{energy})--(\ref{entropy}) the corresponding
initial values $e^0$, $E_\Gamma^0$, and $s^0$ for specific energy,
boundary energy, and entropy, respectively.
Let $E^0 = \io \varrho_0 e^0\dd x$, $S^0 = \io \varrho_0 s^0 \dd x$
denote the total initial energy and entropy, respectively.
{}From the energy end entropy balance equations (\ref{princ1}), (\ref{princ3}), we derive
the following crucial (formal for the moment) balance equation for the ``extended'' energy
$\varrho_0 (e - \theta_\Gamma s)$:
\begin{eqnarray}\label{crucial}
&&\hspace{-16mm}\io\left(c\theta + \frac{\lambda}{2}(U - \alpha(1-\chi))^2
+ \beta \theta_c U + L\chi \right)(x,t)\dd x
+ \frac{K_{\Gamma}}{2} \left(U_{\Omega}(t) + \frac{p_0}{K_{\Gamma}}\right)^2 \\ \nonumber
&&+\, \theta_\Gamma \int_0^t\io \left(\frac{\kappa |\nabla\theta|^2}{\theta^2}
+ \frac{\gamma}{\theta} \chi_t^2 + \frac{\nu}{\theta}U_t^2  \right)(x,\tau)\dd x \dd\tau\\ \nonumber
&&+\, \int_0^t \ipo\frac{h(x)}{\theta}(\theta_\Gamma - \theta)^2(x,\tau)\dd s(x) \dd\tau\\ \nonumber
&=& E^0 + E_\Gamma^0 - \theta_\Gamma S^0 + \theta_\Gamma\io\left(c\log\left(\frac{\theta}{\theta_c}\right)  +
\frac{L}{\theta_c}\chi + \beta U \right)(x,t)\dd x\,.
\end{eqnarray}
We have $\log(\theta/\theta_c) = \log(\theta/2\theta_\Gamma) - \log(\theta_c/2\theta_\Gamma)
\le (\theta/2\theta_\Gamma) - 1 - \log(\theta_c/2\theta_\Gamma)$, hence there exists
a constant $C>0$ independent of $t$ such that for all $t>0$ we have
\begin{eqnarray}\label{esti1}
&&\hspace{-16mm}\io\left(\theta + U^2\right)(x,t)\dd x
+ \int_0^t\io \left(\frac{|\nabla\theta|^2}{\theta^2}
+ \frac{\chi_t^2}{\theta}  + \frac{U_t^2}{\theta}  \right)(x,\tau)\dd x \dd\tau\\ \nonumber
&&+\, \int_0^t \ipo\frac{h(x)}{\theta}(\theta_\Gamma - \theta)^2(x,\tau)\dd s(x) \dd\tau \ \le \ C\,.
\end{eqnarray}

\section{Equilibria}\label{equi}

It follows from (\ref{bctheta}) and (\ref{sys2}) that the only possible equilibrium
temperature is $\theta = \theta_\Gamma$, and the equilibrium configurations
$U_\infty, \chi_\infty$ for $U, \chi$ satisfy for a.e. $x \in \Omega$ the equations
\begin{eqnarray}\label{eqsys1}
\lambda U_\infty(x) - \alpha\lambda (1-\chi_\infty(x)) &=&
\beta (\theta_\Gamma-\theta_c)  -  p_0 - K_{\Gamma}\io U_\infty(x') \dd x'\,,\\ [2mm]\label{eqsys3}
-\lambda U_\infty(x) + \alpha\lambda (1-\chi_\infty(x)) &\in&  \frac{L}{\alpha}\left(1 - \frac{\theta_\Gamma}{\theta_c}\right)
+ \partial I(\chi_\infty(x))\,,
\end{eqnarray}
as a consequence of (\ref{sys1}), (\ref{sys3}), hence
\begin{equation}\label{eqsys2}
\frac{L}{\alpha}\left(\frac{\theta_\Gamma}{\theta_c}-1\right)
- \beta (\theta_\Gamma-\theta_c) + p_0
+ K_{\Gamma}\io U_\infty(x') \dd x'\ \in\ \partial I(\chi_\infty(x))\quad \mbox{a.e.}
\end{equation}
The equilibrium pressure $P_\infty$ is given by (\ref{bcu2}), that is,
\begin{equation}\label{pres}
P_\infty = p_0 + K_{\Gamma}\io U_\infty(x') \dd x'\,.
\end{equation}
Integrating Eq.~(\ref{eqsys1}) over $\Omega$ yields
\begin{equation}\label{eqsys4}
(\lambda + K_{\Gamma} |\Omega|)\io U_\infty(x') \dd x'
= |\Omega|(\beta (\theta_\Gamma-\theta_c) - p_0)
+ \alpha\lambda \io (1-\chi_\infty(x')) \dd x'\,.
\end{equation}
Hence, a necessary and sufficient condition for $\chi_\infty(x)$ to be an equilibrium
phase distribution reads
\begin{equation}\label{eqsys5}
\frac{L}{\alpha\lambda}\left(\frac{\theta_\Gamma}{\theta_c}-1\right)
-\frac{\beta (\theta_\Gamma-\theta_c)  -  p_0}{\lambda + K_{\Gamma} |\Omega|}
+\frac{\alpha K_{\Gamma}}{\lambda + K_{\Gamma} |\Omega|}\io (1-\chi_\infty(x')) \dd x'
\ \in\ \partial I(\chi_\infty(x))\quad \mbox{a.e.}
\end{equation}
Let us introduce a positive dimensionless parameter
\begin{equation}\label{1d17}
{d} := \frac{\alpha^2\lambda K_{\Gamma} |\Omega|}{L (\lambda + K_{\Gamma} |\Omega|)}\,.
\end{equation}
Assume first that $\beta/(\lambda + K_{\Gamma} |\Omega|)$ and
$p_0/(\lambda + K_{\Gamma} |\Omega|)$ are negligible with respect to the other terms.
We then rewrite Eq.~(\ref{eqsys5}) in a simpler form
\begin{equation}\label{1d15}
\frac{\theta_\Gamma}{\theta_c}-1 + \frac{d}{|\Omega|} \io (1-\chi_\infty(x')) \dd x'
\ \in\ \partial I(\chi_\infty(x))\quad \mbox{a.e.}
\end{equation}
We distinguish three cases:

\begin{description}
\item{\fbox{$\theta_\Gamma \ge \theta_c$}}\\[2mm]
Then (\ref{1d15}) can only be satisfied if $\chi_\infty = 1$ a.e.,
hence, by \eqref{eqsys4}, $U_\infty = 0$ a.e.,
and by \eqref{pres}, the pressure $P_\infty$ is in equilibrium with the external
pressure.
We only have the liquid phase in $\Omega$ and the system is stress-free.
\item{\fbox{${d}<1$ and $\theta_\Gamma \le (1-{d})\theta_c$}}\\[2mm]
Then, similarly, (\ref{1d15}) can only be satisfied if $\chi_\infty = 0$ a.e., hence
$$
\io U_\infty(x') \dd x' = \frac{\alpha\lambda |\Omega|}{\lambda + K_{\Gamma} |\Omega|}\,,
\quad U_\infty(x) = \alpha - \frac{\alpha K_{\Gamma} |\Omega|}{\lambda + K_{\Gamma} |\Omega|} =
\frac{\alpha\lambda}{\lambda + K_{\Gamma} |\Omega|}\,.
$$
We only have the solid phase subject to a balance between a positive volume
expansion $U_\infty$ and pressure $P_\infty - p_0 = K_{\Gamma}|\Omega| U_\infty$.

In the limit case $K_{\Gamma} \to 0$ (stress-free
boundary condition, i.e. infinitely soft bottle), we get $P_\infty \to p_0$,
$U_\infty \to \alpha$, ${d} \to 0$. Hence, $\alpha$ measures indeed the
relative volume expansion in the stress-free case.
Similarly, in the limit case $K_{\Gamma} \to \infty$ (rigid bottle), we have
$P_\infty - p_0 \to \alpha \lambda$, $U_\infty \to 0$.
In this case, $\alpha\lambda$ is the pressure difference between inside and outside the bottle.
\item{\fbox{$(1-{d})\theta_c < \theta_\Gamma < \theta_c$}}\\[2mm]
Set ${d}_* = 1 - (\theta_\Gamma/\theta_c) < {d}$. Then every function
$\chi_\infty$ with values in $[0,1]$ satisfying the condition
$(1/|\Omega|)\io (1-\chi_\infty(x'))\dd x' = {d}_*/{d}$ is an equilibrium.
Hence, in this temperature range, we have a large number of possible equilibria.
\end{description}

We thus observe stable undercooled mushy regions in a nonzero temperature range,
and full solidification only takes place if the temperature is below the value
$(1-{d})\theta_c$. Theoretically, we cannot exclude the case ${d} \ge 1$, which
would mean that the solid phase can never be achieved. We show now that in the case
of water and ice, which is relevant for applications, the undercooling coefficient ${d}$
is less than $1$. Approximate values of the physical constants are listed in Table \ref{tablc},
see \cite{ger}.

The maximum of ${d}$ is achieved in a rigid bottle (i.e. $K_{\Gamma} \to \infty$).
By Table \ref{tablc} we have
$\alpha = (V_{ice}-V_{water})/V_{water} = 0.09$, $\lambda \approx 2.25\cdot 10^9\,J/m^3$.
Using Eq.~(\ref{zero}) we obtain
$L = \varrho_0 L_0 \approx 0.33\cdot 10^9\,J/m^3$, hence
${d} = \alpha^2\lambda/L \approx 5.5\%$.
Note that the standard atmospheric pressure is about $10^{5}\, J/m^3$,
while the pressure inside the bottle attains $\alpha\lambda \approx 2\cdot 10^8\,J/m^3$.
This corresponds to a mass of $20$ kilograms pressing by gravity on each square millimeter.

\begin{table}[ht]
\begin{tabular}{|l|c|c|c|}
\hline
{Specific volume of water} & $V_{water} = 1/\varrho_0$ & $10^{-3}$ & $m^{3}/kg$\\
{Specific volume of ice} & $V_{ice}$ & $1.09\cdot 10^{-3}$ & $m^{3}/kg$\\
{Speed of sound} & $v_0 = \sqrt{\lambda/\varrho_0} $   & $ 1.5\cdot 10^3$ & $m/s$\\
{Freezing point} & $\theta_c$ & $273$ & ${K}$\\
{Specific heat} & $c_0$   & $4.2\cdot 10^3$ & $J/(kg\,{K})$\\
{Latent heat} & $L_0$   & $3.3\cdot 10^5$ & $J/kg$\\
{Thermal expansion coefficient} & $\beta/\lambda$ & $\ \, 2.0\cdot 10^{-4}$ & $ {K}^{-1}$\\
\hline
\end{tabular}
\caption{Physical constants for water}
\label{tablc}
\end{table}

In reality, some values of the constants are different in water and in ice
(the specific heat, for instance, is only $2\cdot 10^3\,J/(kg\,{K})$ in the ice).
A phase field model without mechanical effects for this situation was considered in \cite{krs2}.
Also the speed of sound in ice is about the double of the one in water.
We can in principle state the problem with coefficients depending
on $\theta$ and $\chi$ here, too, but this would
lead to serious technical difficulties that we want to avoid here.
Moreover, in water and ice, the thermal expansion coefficient $\beta$ is not constant
and depends strongly on the temperature as well as on the phase.
It may even become negative for temperatures in a right neighborhood of the freezing point.
The values given in Table \ref{tablc} are obtained by a rough linearization
in order to have an idea about the orders of magnitude.

For the coefficient $\beta$ we compute the estimate
$\beta/\varrho_0 = v_0^2\beta/\lambda \approx
450 \,J/(kg\,{K})$, while $L_0/(\alpha\theta_c) = L/(\varrho_0\alpha\theta_c)
\approx 13400 \,J/(kg\,{K})$. Let us define a new constant
$L_\beta := L_0 - \beta\alpha\theta_c/\varrho_0$
as a small (3.3\%) correction to the latent heat $L_0$. Using (\ref{pres}), we may rewrite
(\ref{eqsys2}) as
\begin{equation}\label{eqsys2c}
\frac{\varrho_0 L_\beta}{\alpha}\left(\frac{\theta_\Gamma}{\theta_c}-1\right)
+ P_\infty\ \in\ \partial I(\chi_\infty(x))\quad \mbox{a.e.}
\end{equation}
We now show that (\ref{eqsys2c}) contains
the Clausius-Clapeyron equation,  cf.~\cite[Book~5, Chapter~5]{joo} or \cite[pp.~124--126]{ye}.
The pressure $P_\infty$
is defined as the difference $\delta P$ between the absolute pressure
and the standard pressure. The phase transition takes place at temperature $\theta_\Gamma$
if the right hand side of (\ref{eqsys2c}) vanishes. The temperature difference
is $\delta\theta = \theta_\Gamma - \theta_c$, and we get the Clausius-Clapeyron relation
in the form of Eq.~(288) of \cite{ye}, that is,
\begin{equation}\label{clacla}
\frac{\delta P}{\delta\theta} = -\frac{\varrho_0 L_\beta}{\alpha\theta_c} =
\frac{L_\beta}{\theta_c (V_{water}-V_{ice})}\,.
\end{equation}

For general $\beta \ge 0$ and $p_0$, we have an analogous classification as above.
We introduce further dimensionless quantities
\begin{equation}\label{1d17a}
\tilde\beta = \frac{\alpha\lambda \beta\theta_c}{L (\lambda + K_{\Gamma} |\Omega|)}\,,\quad
\omega = \frac{\alpha\lambda  p_0}{L (\lambda + K_{\Gamma} |\Omega|)}\,.
\end{equation}
The counterpart of (\ref{1d15}) reads
\begin{equation}\label{1d15a}
(1-\tilde\beta)\left(\frac{\theta_\Gamma}{\theta_c}-1\right) + \omega
+ \frac{d}{|\Omega|} \io (1-\chi_\infty(x')) \dd x'
\ \in\ \partial I(\chi_\infty(x))\,.
\end{equation}
Assuming that $\tilde\beta < 1$, we thus observe pure liquid
for $\theta_\Gamma \ge \theta_c(1-\omega/(1-\tilde\beta))$,
while pure solid corresponds to $\theta_\Gamma \le \theta_c(1-(\omega + d)/(1-\tilde\beta))$.
The dimensionless external pressure deviation $\omega$ can be assumed small. However,
$\tilde\beta$ is a material constant, and
the condition $\tilde\beta < 1$ might be restrictive. Again, for water and ice,
the maximal value $\tilde\beta = (L_0-L_\beta)/L_0 \approx 0.033$
corresponding to $K_{\Gamma} = 0$ shows that the influence of thermal expansion
on the undercooling coefficient is negligible.

\section{Existence and uniqueness of solutions}\label{exiu}

We construct the solution of (\ref{sys1})--(\ref{sys3}) by the Banach
contraction argument. The method of proof is independent
of the actual values of the material constants, and we choose for simplicity
\begin{equation}\label{const}
L= 2,\ \ c= \theta_c =\alpha = \beta= \gamma = \kappa = \lambda = \nu = 1\,.
\end{equation}
System (\ref{sys2})--(\ref{sys3}) with boundary condition (\ref{bctheta}) then reads
\begin{eqnarray}\label{nequ1}
\io\theta_t w(x)\dd x + \io \nabla\theta \cdot \nabla w(x) \dd x &=&
\io\left(U_t^2 - \theta U_t -\Big(U +\chi +1\Big)\chi_t\right) w(x)\dd x\qquad \\ \nonumber
&&-\,\ipo h(x)(\theta - \theta_\Gamma) w(x) \dd s(x) \,,\\ [2mm]\label{nequ2}
U_t + U + \chi + K_{\Gamma}U_{\Omega}(t)  &=&  \theta  - p_0\,,\\ [2mm]\label{nequ3}
\chi_t + U +\chi + \partial I(\chi) &\ni& 2\theta - 1\,,
\end{eqnarray}
where (\ref{nequ1}) is to be satisfied for all test functions $w \in W^{1,2}(\Omega)$
and a.e. $t>0$, while (\ref{nequ2})--(\ref{nequ3}) are supposed to hold a.e. in
$\Omega_\infty := \Omega\times (0,\infty)$.

In this section we prove the following existence and uniqueness result.
\begin{theorem}\label{main}
Let $0< \theta_* \le \theta_\Gamma \le \theta^*$ and $p_0 \in \real$
be given constants, and let the data satisfy the conditions
$$
\barr{rcllcccll}
\theta^0 &\in& W^{1,2}(\Omega)\cap L^\infty(\Omega)\,,
\quad &\theta_* &\le& \theta^0(x)&\le& \theta^*\quad &\mbox{a.e.}\,,\\
U^0, \chi^0 &\in& L^\infty(\Omega)\,,\quad &0 &\le& \chi^0(x) &\le& 1 \quad &\mbox{a.e.}
\earr
$$
Then there exists a unique solution $(\theta,U,\chi)$ to (\ref{nequ1})--(\ref{nequ3}),
(\ref{ini2})--(\ref{ini3}), such that
$\theta>0$ a.e., $\chi \in [0,1]$ a.e., 
$U, U_t,\chi_t, \theta, 1/\theta \in L^\infty(\Omega_\infty)$,
$\theta_t, \Delta\theta \in L^2(\Omega_\infty)$,
and $\nabla\theta \in L^\infty(0,T;L^2(\Omega)) \cap L^2(\Omega_\infty)$.
\end{theorem}

\begin{remark}\label{bc}
\rm For existence and uniqueness alone, we might allow the external temperature
$\theta_\Gamma$ to depend on $x$ and $t$, and assume only that it belongs to the space
$W^{1,2}_{{\rm loc}} (0,\infty;L^2(\partial\Omega))
\cap L^\infty_{{\rm loc}}(\partial\Omega\times(0,\infty))$.
For the global bounds, the assumption that $\theta_\Gamma$ be constant plays a
substantial role.
\end{remark}

The proof of Theorem \ref{main} will be carried out in the following subsections.
Notice first that the term $U_t^2 - \theta U_t - (U +\chi +1)\chi_t$ on the right hand
side of (\ref{nequ1}) can be rewritten alternatively, using \eqref{nequ3} and \eqref{nequ2}, as
\begin{eqnarray}\label{rhs}
U_t^2 - \theta U_t - (U +\chi +1)\chi_t &=& U_t^2 - \theta U_t + \chi_t^2 - 2 \theta\chi_t
\\[2mm]\nonumber
&=& - (\chi + U + p_0 + K_{\Gamma}U_{\Omega}) U_t - (U +\chi +1)\chi_t\,,
\end{eqnarray}
We now fix some constant $R>0$ and construct the solution for the truncated system
\begin{eqnarray}\label{trequ1}
\io\theta_t w(x)\dd x + \io \nabla\theta \cdot \nabla w(x) \dd x &=&
\io\left(U_t^2 +\chi_t^2 -  Q_R(\theta) (U_t+ 2 \chi_t)\right) w(x)\dd x\qquad \\ \nonumber
&&\hspace{-15mm}-\,\ipo h(x)(\theta - \theta_\Gamma) w(x)\dd s(x)\quad \forall w \in W^{1,2}(\Omega)\,,
\\ [2mm]\label{trequ2}
U_t + U +  \chi + K_{\Gamma}U_\Omega(t) &=& Q_R(\theta) - p_0\,,\\ [2mm]\label{trequ3}
 \chi_t + U +\chi + \partial I(\chi) &\ni& 2 Q_R(\theta) - 1
\end{eqnarray}
first in a bounded domain $\Omega_T := \Omega\times (0,T)$ for any given $T>0$,
where $Q_R$ is the cutoff function $Q_R(z) = \min\{z^+, R\}$. We then derive upper and lower
bounds for $\theta$ independent of $R$ and $T$, so that the local solution of 
(\ref{trequ1})--(\ref{trequ3}) is also a global solution of (\ref{nequ1})--(\ref{nequ3})
if $R$ is sufficiently large.

\subsection{A gradient flow}\label{grad}

In a separable Hilbert space $H$ with norm $|\cdot|$,
consider a gradient flow
\begin{equation}\label{eg1}
\dot v(t) + \partial \psi(v(t)) \ \ni \ f(t)\,, \quad v(0) = v^0\,,
\end{equation}
where $\psi: H \to [0,\infty]$ is a proper convex lower semicontinuous functional
such that $\lim_{|v| \to \infty} \psi(v) = +\infty$,
$\partial \psi$ is its subdifferential, and $v^0 \in {\rm Dom\,}\psi$, $f \in L^2(0,\infty;H)$
are given. A classical existence and uniqueness result in \cite[Th\'eor\`eme 3.6]{B}
states that for every $T>0$ there exists a unique solution $v \in C([0,T];H)$ to (\ref{eg1})
such that $\dot v \in L^2(0,T;H)$, and
$$
\left(\int_0^T |\dot v(\tau)|^2\dd \tau\right)^{1/2} \le \psi(v^0) +
\left(\int_0^T |f(\tau)|^2\dd \tau\right)^{1/2}\,.
$$
 We prove here the following Lemma.

\begin{lemma}\label{lg1}
Let $f, \dot f$ belong to $L^2(0,\infty;H)$. Then $\lim_{t \to \infty} \dot v(t) = 0$.
\end{lemma}

\bpf{Proof}\
For each $h>0$ and a.e. $t>0$ we have
$$
\frac 12 \frac{\dd}{\dd t} \left|\frac{v(t+h) - v(t)}{h}\right|^2
\ \le \ \left|\frac{f(t+h) - f(t)}{h}\right|\,\left|\frac{v(t+h) - v(t)}{h}\right|\,,
$$
hence
\begin{equation}\label{eg2}
|\dot v(t)|^2 - |\dot v(s)|^2
\ \le \ 2\int_s^t |\dot f(\tau)| |\dot v(\tau)|\dd \tau
\end{equation}
for almost all $0< s < t$. Hence, the function $t \mapsto
2\int_0^t |\dot f(\tau)| |\dot v(\tau)|\dd \tau
- |\dot v(t)|^2$ is almost everywhere equal to a nondecreasing function in $(0,\infty)$.
We are thus in the situation of \cite[Proposition 5.2]{ksz}, which
gives the desired statement.
\epf

We apply the above result to the case $H=L^2(\Omega)\times L^2(\Omega)$, and
\begin{eqnarray}\label{eg4}
v &=& \vect{U}{\chi}\,,\\[2mm]\label{eg5}
 \quad \psi(v) &=& \io \left( \frac12 (U+\chi - 1)^2
+ (U + 2\chi)(1 - \theta_\Gamma) + I(\chi)\right)\dd x\\ \nonumber
&&+\, \frac{K_{\Gamma}}{2}\left(\io U \dd x + \frac{p_0}{K_{\Gamma}}\right)^2 + C_\psi\,,\\[2mm]\label{eg6}
f &=& \vect{Q_R(\hat\theta) - \theta_\Gamma}{2(Q_R(\hat\theta) - \theta_\Gamma)}\,,
\end{eqnarray}
where $C_\psi$ is a suitable constant such that $\psi(v)\ge 0$ for all $v$, and
$\hat \theta$ is a given function. The initial condition $v^0$
is given by (\ref{ini1}), (\ref{ini3}). We have
\begin{equation}\label{eg7}
\vect{\eta}{\zeta} \in \partial\psi(v) \ \Longleftrightarrow \
\left\{
\begin{array}{rcl}
\eta &=& U+\chi - \theta_\Gamma + K_{\Gamma}\io U\dd x + p_0\,,\\
\zeta &\in& U+\chi + 1  - 2\theta_\Gamma + \partial I(\chi)\,,
\end{array}
\right.
\end{equation}
and we see that Eqs.~(\ref{trequ2})--(\ref{trequ3}) with $\theta$ replaced by $\hat\theta$
can be equivalently written as a gradient flow (\ref{eg1}), (\ref{eg4})--(\ref{eg6}).
For its solutions, we prove the following result.

\begin{proposition}\label{pg1}
Let the hypotheses of Theorem \ref{main} hold, and let a function
$\hat\theta \in L^2_{{\rm loc}}(0,\infty; L^2(\Omega))$ be given.
Let $(U,\chi)$ be the solution of
(\ref{eg1}), (\ref{eg4})--(\ref{eg6}). Then there exists a constant $C_0$, independent
of $x, t$ and $R$, such that a.e. in $\Omega_\infty$ we have
\begin{equation}\label{eg8bis}
|U(x,t)| + |U_t(x,t)| + |\chi_t(x,t)| \le C_0(1+R)\,.
\end{equation}
Let furthermore $\hat\theta_1, \hat\theta_2 \in L^2_{{\rm loc}}(0,\infty; L^2(\Omega))$
be two functions, and let $(U_1,\chi_1), (U_2,\chi_2)$ be the corresponding solutions of
(\ref{eg1}), (\ref{eg4})--(\ref{eg6}). Then the differences
$\hat\theta_d = \hat\theta_1 - \hat\theta_2$, $U_d = U_1-U_2$, $\chi_d =\chi_1-\chi_2$
satisfy for every $t \ge 0$ and a.e. $x \in \Omega$ the inequality
\begin{equation}\label{eg9}
\int_0^t (|(U_d)_t| + |(\chi_d)_t|)(x,\tau) \dd \tau \le C_0(1+t)
\int_0^t\left(|\hat\theta_d(x,\tau)| + t |\hat\theta_d(\tau)|_2\right)\dd \tau\,,
\end{equation}
where the symbol $|\cdot|_2$ stands for the norm in $L^2(\Omega)$.
\end{proposition}

In what follows,
we denote by $C_1, C_2, \dots$ any constant independent of $x, t$ and~$R$.

\bpf{Proof}
Put $X_{\Omega}(t) = \io (1-\chi(x',t))\dd x'$.
Integrating (\ref{trequ2}) with $\theta$ replaced by $\hat\theta$ over $\Omega$ yields
$$
\dot U_{\Omega} + (1+K_{\Gamma} |\Omega|) U_{\Omega} = X_{\Omega}
+ \io (Q_R(\hat\theta)-1) \dd x - p_0 |\Omega|\quad \mae
$$
Since $\chi$ attains values in $[0,1]$, we easily obtain
$|\dot U_{\Omega}| + |U_{\Omega}| \le C_1(1+R)\ \mae$ Equation (\ref{trequ2}) now
has a right hand side bounded by a multiple of $1+R$, hence $|U_t| + |U| \le C_2(1+R)\ \mae$
To obtain the same bound
for $|\chi_t|$, it suffices to multiply (\ref{trequ3}) by $\chi_t$.
This completes the proof of \eqref{eg8bis}.

To prove (\ref{eg9}), we rewrite (\ref{eg1}), (\ref{eg4})--(\ref{eg6}) as two scalar gradient
flows
\begin{eqnarray}\label{eg10}
U_t + \partial \psi_1(U) &=& a\,,\\[1mm]\label{eg11}
\chi_t + \partial \psi_2(\chi) &\ni& b\,,
\end{eqnarray}
where $\psi_1(U) = \frac12 U^2$, $\psi_2 = \frac12 \chi^2 + I(\chi)$,
$a = Q_R(\hat\theta) -\chi - p_0 - K_{\Gamma} U_{\Omega}$, $b = 2 Q_R(\hat\theta)-1-U$.
Consider now two different inputs. As above, we denote the differences $\{\}_1-\{\}_2$
by $\{\}_d$ for all symbols $\{\}$.
By \cite[Theorem 1.12]{hys}, we have for all $t>0$ and a.e. $x \in \Omega$ that
 \begin{equation}\label{eg12}
\int_0^t (|(U_d)_t| + |(\chi_d)_t|)(x,\tau) \dd \tau \le 2\int_0^t
(|a_d| + |b_d|)(x,\tau) \dd \tau\,.
 \end{equation}
We multiply the difference of (\ref{eg10}) by $U_d$, the difference of (\ref{eg11})
by $\chi_d$, and sum them up to obtain that
\begin{equation}\label{eg13}
(U_d)_t U_d + (\chi_d)_t\chi_d + (U_d +\chi_d)^2 + K_{\Gamma} U_{\Omega d}\,U_d
\le |\hat\theta_d|(|U_d| + 2 |\chi_d|) \quad\mae
\end{equation}
We first integrate (\ref{eg13}) over $\Omega$. Using the symbol $|\cdot|_2$ for the norm
in $L^2(\Omega)$, we get for a.e. $t >0$ that
\begin{equation}\label{eg14}
\frac12 \frac{\dd}{\dd t} \left(|U_d|_2^2 + |\chi_d|_2^2\right) + K_{\Gamma}U_{\Omega d}^2 \le
|\hat\theta_d|_2(|U_d|_2 +2 |\chi_d|_2) \le
\sqrt{5}|\hat\theta_d|_2 \left(|U_d|_2^2 +|\chi_d|_2^2\right)^{1/2}\,.
\end{equation}
Hence, $\frac{\dd}{\dd t} (|U_d|_2^2 + |\chi_d|_2^2)^{1/2} \le \sqrt{5}|\hat\theta_d|_2$
a.e., and integrating over $t$, we find that
 \begin{equation}\label{eg15}
\left(|U_d|_2^2 +|\chi_d|_2^2\right)^{1/2}(t)
\le \sqrt{5}\int_0^t |\hat\theta_d(\tau)|_2\dd\tau\,.
 \end{equation}
This implies in particular that
\begin{equation}\label{eg16}
|U_{\Omega d}(t)| \le \sqrt{5|\Omega|}\int_0^t |\hat\theta_d(\tau)|_2\dd\tau\,.
\end{equation}
Using again (\ref{eg13}), we find for a.e. $(x,t)\in \Omega_\infty$ the inequality
\begin{equation}\label{eg17}
\frac12 \frac{\partial}{\partial t} \left(|U_d|^2 + |\chi_d|^2\right)(x,t)  \le
\left(K_{\Gamma}|U_{\Omega d}(t)|+ |\hat\theta_d(x,t)|\right)(|U_d| +2 |\chi_d|)(x,t)\,.
\end{equation}
This is for almost all $x\in\Omega$ an inequality of the form
$(\dd/\dd t)(Y^2(t)) \le 2 c(t)Y(t)$, $Y(0) = 0$, with $Y = (|U_d|^2 + |\chi_d|^2)^{1/2}$,
which implies $Y(t) \le \int_0^t c(\tau)\dd\tau$ for all $t>0$.  Hence,
\begin{equation}\label{eg18}
\left(|U_d|^2 + |\chi_d|^2\right)^{1/2} (x,t)  \le
C_1 \int_0^t \left(|\hat\theta_d(x,\tau)| + t |\hat\theta_d(\tau)|_2\right)\dd\tau
\quad \mbox{a.e.}
\end{equation}
This enables us to estimate the right hand side of (\ref{eg12}) and obtain the bound
 \begin{equation}\label{eg19}
\int_0^t (|(U_d)_t| + |(\chi_d)_t|)(x,\tau)\dd \tau \le C_2\int_0^t
\left((1+t)|\hat\theta_d(x,\tau)| + t(1+t) |\hat\theta_d(\tau)|_2\right)\dd \tau
 \end{equation}
for a.e. $x\in \Omega$ and all $t \ge 0$. This completes the proof.
\epf

\subsection{Existence of solutions for the truncated problem}\label{exitrunc}

We construct the solution of (\ref{trequ1})--(\ref{trequ3}) for every $R>0$ by the Banach
contraction argument on a fixed time interval $(0,T)$.

\begin{lemma}\label{le1}
Let the hypotheses of Theorem \ref{main} hold, and let $T>0$ and $R>0$ be given.
Then there exists
a unique solution $(\theta,U,\chi)$ to (\ref{trequ1})--(\ref{trequ3}),
(\ref{ini2})--(\ref{ini3}), such that $U\in W^{1,\infty}(\Omega_T)$, $\theta>0$ a.e.,
$\chi_t, \theta, 1/\theta \in L^\infty(\Omega_T)$, $\theta_t, \Delta\theta
\in L^2(\Omega_T)$, and $\nabla\theta \in L^\infty(0,T;L^2(\Omega))$.
\end{lemma}

\bpf{Proof}\ Let $\hat\theta \in L^2(\Omega_T)$ be a given function,
and consider the system
\begin{eqnarray}\label{appqu1}
\io\theta_t w(x)\dd x + \io \nabla\theta \cdot \nabla w(x) \dd x &=&
\io (U_t^2 +\chi_t^2 -  Q_R(\hat \theta) (U_t+2\chi_t)) w(x)\dd x\qquad\quad \\ \nonumber
&&\hspace{-15mm}-\,\ipo h(x)(\theta - \theta_\Gamma) w(x) \dd s(x) \quad \forall w \in W^{1,2}(\Omega)\,,
\\ [2mm]\label{appqu2}
U_t + U +\chi + K_{\Gamma}U_\Omega(t) &=& Q_R(\hat\theta)- p_0\,,\\[2mm]\label{appqu3}
 \chi_t +  U +\chi + \partial I(\chi) &\ni& 2 Q_R(\hat\theta) - 1\,.
\end{eqnarray}
Equations (\ref{appqu2})--(\ref{appqu3}) are solved as a gradient flow problem from
Subsection \ref{grad}, while (\ref{appqu1}) is a simple linear parabolic equation for $\theta$.
Testing (\ref{appqu1}) by $\theta_t$, we obtain by Proposition \ref{pg1} that
\begin{eqnarray}\label{eg19a}
&&\int_0^T\io \theta_t^2\dd x \dd t + \supess_{t\in (0,T)}\left(
\io |\nabla\theta|^2 \dd x + \ipo h(x)(\theta - \theta_\Gamma)^2 \dd s(x)\right)
\\[2mm]\nonumber
&& \qquad \le T |\Omega| \left(C_0(1+R)(2 C_0(1+R) + 3R)\right)^2 =: M_R\,.
\end{eqnarray}
Hence, we can define the mapping that with $\hat\theta$ associates the solution $\theta$
of (\ref{appqu1})--(\ref{appqu3}) with initial conditions (\ref{ini2})--(\ref{ini3}).
We now show that it is a contraction on the set
\begin{equation}\label{eg8}
\Xi_{T,R} := \{\hat\theta \in L^2(\Omega_T): \mbox{conditions (\ref{eg8a})--(\ref{eg8d})
hold}\}\,,
\end{equation}
where
\begin{eqnarray}\label{eg8a}
&&\hspace{-14mm}\hat\theta_t \in L^2(\Omega_T)\,;\\[2mm]\label{eg8b}
&&\hspace{-14mm}\nabla\hat\theta \in L^\infty(0,T;L^2(\Omega))\,;\\[2mm]\label{eg8c}
&&\hspace{-14mm}\int_0^T\io \hat\theta_t^2 \dd x \dd t + \supess_{t\in (0,T)}\left(
\io |\nabla\hat\theta|^2 \dd x + \ipo h(x)(\hat\theta - \theta_\Gamma)^2 \dd s(x)\right)
\le M_R\,;\\[2mm]\label{eg8d}
&&\hspace{-14mm}\hat\theta(x,0) = \theta^0(x)\ \mae
\end{eqnarray}
Let $\hat\theta_1, \hat\theta_2$
be two functions in $\Xi_{T,R}$, and let $(\theta_1, U_1, \chi_1)$, $(\theta_2, U_2, \chi_2)$,
be the corresponding solutions to (\ref{appqu1})--(\ref{appqu3}) with the same initial
conditions $\theta^0, U^0,\chi^0$. We see from (\ref{eg19a}) that 
$\theta_1, \theta_2$ belong to $\Xi_{T,R}$.
Integrating Eq.~(\ref{appqu1}) for $\theta_1$ and $\theta_2$ with respect to time and testing
their difference by $w = \theta_d := \theta_1-\theta_2$, we obtain, using Proposition~\ref{pg1}, that
\begin{eqnarray}\nonumber
&& \hspace{-15mm}\io\theta_d^2(x,t)\dd x +
\frac{\dd}{\dd t}\left(\io\left|\nabla\int_0^t \theta_d(x,\tau)\dd\tau\right|^2\dd x
+ \ipo h(x) \left|\int_0^t \theta_d(x,\tau)\dd\tau\right|^2\dd s(x)\right)\\[2mm]  \label{lipthe2}
 &\le& C_3(1+R) \io \left(\int_0^t (|(U_d)_t| + |(\chi_d)_t| + |\hat\theta_d|)
(x,\tau) \dd \tau\right) \theta_d(x,t) \dd x\quad \mae
\end{eqnarray}
{}From (\ref{eg9}) and Minkowski's inequality, it follows that
\begin{eqnarray*}
\left|\int_0^t (|(U_d)_t| + |(\chi_d)_t|)(\cdot,\tau) \dd \tau\right|_2
&\le& C_4 (1+t)^2\int_0^t |\hat\theta_d(\tau)|_2 \dd \tau\\[2mm]
&\le& C_4 (1+t)^2\left(t\int_0^t |\hat\theta_d(\tau)|_2^2 \dd \tau\right)^{1/2}.
\end{eqnarray*}
By Young's inequality, we rewrite (\ref{lipthe2}) as
\begin{eqnarray}\nonumber
&& \hspace{-15mm}\io\theta_d^2(x,t)\dd x +
\frac{\dd}{\dd t}\left(\io\left|\nabla\int_0^t \theta_d(x,\tau)\dd\tau\right|^2\dd x
+ \ipo h(x) \left|\int_0^t \theta_d(x,\tau)\dd\tau\right|^2\dd s(x)\right)\\[2mm]  \label{lipthe4}
 &\le& C_5(1+R^2)(1+t)^5\int_0^t |\hat\theta_d(\tau)|_2^2 \dd \tau \quad \mae
\end{eqnarray}
Set $\Theta^{2} (t) = \int_0^t|\theta_d(\tau)|_2^2 \dd\tau$,
$\hat\Theta^{2} (t) = \int_0^t|\hat\theta_d(\tau)|_2^2 \dd\tau$.
Integrating (\ref{lipthe4}) with respect to time, we obtain
\begin{equation}\label{lipthe3}
\Theta^{2}(t) \le  C_5(1+R^2) \int_0^t (1+\tau)^5 \hat\Theta^{2}(\tau)\dd\tau\,.
\end{equation}
We set $C_R:= (C_5(1+R^2)/{6})$ and introduce in $L^\infty(0,T)$ the norm
$$
\|w\|_C := \sup_{\tau \in [0,T]} \expe^{- C_R(1+\tau)^{6}}|w(\tau)|\,.
$$
Then $\|\Theta\|_C^{2} \le \frac12\|\hat\Theta\|_C^{2}$, and
hence the mapping $\hat\theta\mapsto \theta$
is a contraction in $L^2(\Omega_T)$ with respect to the norm induced by $\|\cdot\|_C$.
The set $\Xi_{T,R}$ is a closed subset of $L^2(\Omega_T)$. This implies the existence of a 
fixed point $\theta \in \Xi_{T,R}$, which is indeed a solution to (\ref{trequ1})--(\ref{trequ3}).
The positive upper and lower bounds for $\theta$ follow from the maximum principle.
Indeed, the right hand side (\ref{rhs}) of (\ref{trequ1}) is bounded from above by $C_6(1+R)^2$
and from below by $-\frac12 (\theta^+)^2$. Let us define the functions
$$
\theta^\sharp(t) = \theta^* + C_6(1+R)^2t\,,\quad \theta^\flat(t) = \frac{2\theta_*}{2+\theta_* t}\,.
$$
For every nonnegative test function $w$ and a.e. $t\in (0,T)$ we have
\begin{eqnarray}\label{max1}
&&\hspace{-14mm}\io\theta_t w(x)\dd x + \io \nabla\theta \cdot \nabla w(x) \dd x
+\ipo h(x)(\theta - \theta_\Gamma) w(x) \dd s(x)\\[1mm]\nonumber
 &\le&
C_6(1+R)^2\io w(x)\dd x\,,\qquad \\[2mm] \label{max2}
&&\hspace{-14mm}\io\theta_t w(x)\dd x + \io \nabla\theta \cdot \nabla w(x) \dd x
+\ipo h(x)(\theta - \theta_\Gamma) w(x) \dd s(x)\\[1mm]\nonumber
&\ge&
-\frac12\io (\theta^+)^2 w(x)\dd x\,,\qquad \\[2mm] \label{max3}
&&\hspace{-14mm}\io\theta^\sharp_t w(x)\dd x + \io \nabla\theta^\sharp \cdot \nabla w(x) \dd x
+\ipo h(x)(\theta^\sharp - \theta_\Gamma) w(x) \dd s(x)\\[1mm]\nonumber
&\ge&
C_6(1+R)^2\io  w(x) \dd x\,,\qquad \\[2mm] \label{max4}
&&\hspace{-14mm}\io\theta^\flat_t w(x)\dd x + \io \nabla\theta^\flat \cdot \nabla w(x) \dd x
+\ipo h(x)(\theta^\flat - \theta_\Gamma) w(x) \dd s(x)\\[1mm]\nonumber
 &\le&
-\frac12\io (\theta^\flat)^2 w(x)\dd x\,.
\end{eqnarray}
We now subtract (\ref{max3}) from (\ref{max1}) and test by $w = (\theta-\theta^\sharp)^+$,
which yields the pointwise bound $\theta(x,t) \le \theta^\sharp(t)$. Similarly,
we subtract (\ref{max2}) from (\ref{max4}) and test by $w = (\theta^\flat-\theta)^+$.
We thus have the inequalities
\begin{equation}\label{elinf}
\theta^\flat(t)\le\theta(x,t) \le \theta^\sharp(t) \quad \mbox{a.e.},
\end{equation}
which complete the proof of Lemma \ref{le1}.
\epf

\subsection{Proof of Theorem \ref{main}}\label{concluwellp}

The unique solution $(\theta,U,\chi)$ to (\ref{trequ1})--(\ref{trequ3}),
(\ref{ini2})--(\ref{ini3}) exists globally in the whole domain $\Omega_{\infty}$.
We now derive uniform bounds independent of $t$ and $R$.
Take first for instance any $R> 2\theta^*$. By (\ref{elinf}), we know that
the solution component $\theta$ of (\ref{trequ1})--(\ref{trequ3}) remains
smaller than $R$ in a nondegenerate interval $(0,T)$ with $T > \theta^*/(C_6(1+R)^2)$.
Let $(0,T_0)$ be the maximal interval in which $\theta$ is bounded by $R$.
Then, in $(0,T_0)$, the solution given by Lemma~\ref{le1} is also a solution of the original problem
(\ref{nequ1})--(\ref{nequ3}).
Moreover, due to estimate \eqref{esti1}, we know that $\theta$ admits a bound
in $L^\infty(0,T_0;L^1(\Omega))$ independent of $R$. In order to prove that
$T_0 = +\infty$ if $R$ is sufficiently large, we need the following
variant of the Moser iteration lemma.

\begin{proposition}\label{moser}
Let $\Omega \subset \real^N$ be a bounded domain with Lipschitzian boundary.
Given nonnegative functions $h \in L^1(\partial\Omega)$
and $r \in L^\infty(0,\infty; L^{q} (\Omega))$ with a fixed $q > N/2$,
$|r|_{L^\infty(0,\infty; L^q(\Omega))} =: r^*$,
an initial condition $v^0\in L^\infty(\Omega)$, and a boundary datum
$v_\Gamma \in L^\infty(\partial\Omega\times (0,\infty))$, consider the problem
\begin{eqnarray}\label{moser1}
v_t -\Delta v + v &=& r(x,t)\,{\mathcal H}[v]\qquad\hbox{a.e.~in }
\Omega \times (0,\infty)\,,\\
\label{moser1a}
\nabla v\cdot \bfn &=& - h(x)\,\left(f(x,t, v(x,t)) - v_\Gamma(x,t)\right)
\quad\hbox{a.e.~on }\partial\Omega \times (0,\infty)\,,\\
\label{moser2}
v(x,0)&=&v^0\qquad \hbox{a.e.~in }\Omega\,,
\end{eqnarray}
under the assumption that there exist positive constants
$m, H_0, C_f, V, V_\Gamma, E_0$
such that the following holds:
\begin{itemize}
\item[{\rm (i)}] The mapping ${\mathcal H}\,:\,L^\infty_{{\rm loc}}
(\Omega \times (0,\infty))\to
L^\infty_{\rm loc}(\Omega \times (0,\infty))$ satisfies for every
$v\in L^\infty_{\rm loc}(\Omega \times (0,\infty))$ and a.e.
$(x,t) \in \Omega \times (0,\infty)$ the inequality
\[
v(x,t)\,{\mathcal H}[v](x,t) \leq H_0|v(x,t)| \left(1+ |v(x,t)| + \int_0^t
\xi(t-\tau)|v(x,\tau)| \dd\tau\right),
\]
where $\xi \in W^{1,1}(0,\infty)$ is a given nonnegative function such that
\begin{equation}\label{me1a}
\dot \xi(t) \le -\xi(0)\,\xi(t) \quad \mbox{a.e.}
\end{equation}
\item[{\rm (ii)}] $f$ is a Carath\'eodory function on $\Omega \times (0,\infty) \times \real$
such that $f(x,t,v)\,v\geq C_f\,v^2$ a.e. for all $v\in \real$.
\item[{\rm (iii)}] $|v^0(x)|\leq V$ a.e.~in $\Omega$.
\item[{\rm (iv)}] $|v_\Gamma(x,t)|\leq V_\Gamma$ a.e.~on
$\partial\Omega \times (0,\infty)$.
\item[{\rm (v)}] System (\ref{moser1})--(\ref{moser2}) admits a solution\\
$v \in W^{1,2}_{{\rm loc}}(0,\infty;
(W^{1,2})'(\Omega))\cap L^2_{{\rm loc}}(0,\infty;W^{1,2}(\Omega))
\cap L^\infty_{\rm loc}(\Omega \times (0,\infty))$\\
satisfying the estimate
\[
\io|v(x,t)|\,dx\leq E_0 \quad\hbox{a.e.~in }(0,\infty)\,.
\]
\end{itemize}
Then there exists a positive constant $C^*$ depending only on
$|h|_{L^1(\partial\Omega)}$, $C_f$, $H_0$ such that
\begin{equation}\label{moseresti}
|v(t)|_{L^\infty(\Omega)}\leq C^*\max\left\{1, V, V_\Gamma, E_0\right\}\quad
\hbox{for a.e.~}t>0.
\end{equation}
\end{proposition}

\begin{remark}\label{remxi}
{\rm As a consequence of (\ref{me1a}), we have $\xi(t) \le \xi(0)\expe^{-\xi(0)t}$
for all $t\ge 0$, hence $\int_0^\infty \xi(t)\dd t \le 1$. As a typical function
satisfying (\ref{me1a}), let us mention for example
\begin{equation}\label{me1b}
\xi(t) = \frac{m_1}{\sum_{k=1}^n r_k}\sum_{k=1}^n r_k \expe^{-m_k t}
\end{equation}
with any $0< m_1\le  \dots \le m_n$ and $r_k > 0$, $k=1, \dots, n$.}
\end{remark}

We split the proof of Proposition \ref{moser} into several steps.

\begin{lemma}\label{ml1}
Let $\xi$ be as in Proposition \ref{moser}, let
$x \in W^{1,1}_{{\rm loc}}(0,\infty)$ and $y \in L^{1}_{{\rm loc}}(0,\infty)$
be nonnegative functions, and let $a>0$, $C>0$, $\delta \in (0,1)$
be given constants. Set $\mu = \min\{a, \xi(0)(1-\delta)\}$, and assume
that for a.e. $t>0$ we have
\begin{equation}\label{me1}
\dot x(t) + a x(t) + y(t) \le C + \delta
 \int_0^t \xi(t-\tau)\, y(\tau) \dd\tau\,.
\end{equation}
Then $x(t) \le \max\{x(0), C/\mu\}$ for all $t>0$.
\end{lemma}

\bpf{Proof}
Set $z(t) = \int_0^t \xi(t-\tau) y(\tau) \dd\tau$. Then
$(1/\xi(0))\dot z(t) + z(t) \le y(t)$ a.e., hence
$$
\dot x(t) + a x(t) + \frac{1}{\xi(0)} \dot z(t) + (1-\delta) z(t) \le C\quad \mae
$$
With $\mu$ as above, we have
$$
\left(\dot x(t) + \frac{1}{\xi(0)}\dot z(t)\right) +
\mu\left(x(t) + \frac{1}{\xi(0)} z(t)\right) \le C\quad \mae\,,
$$
which yields
$$
x(t) + \frac{1}{\xi(0)} z(t) \le \max\left\{
x(0) + \frac{1}{\xi(0)} z(0), \frac{C}{\mu}\right\}\,,
$$
and the desired inequality follows easily.
\epf

\begin{lemma}\label{ml2}
Let $\mathcal{H}$ be as in Proposition \ref{moser}, and let $|\cdot|_p$ denote
the norm in $L^p(\Omega)$ for $1 \le p \le \infty$. Let
$v\in L^\infty_{\rm loc}(\Omega_\infty)$ and $p,s \ge 1$ be arbitrary.
For $(x,t) \in \Omega_\infty$ set $h_p = \mathcal{H}[v]\, v\, |v|^{p-2}$.
Then, for all $t>0$ we have
$$
|h_p(t)|_s \ \le \ H_0\left(\frac{1}{p} + 3 |v(t)|_{ps}^p + \frac{1}{p}
\int_0^t \xi(t-\tau)|v(\tau)|_{ps}^p \dd\tau\right),
$$
where $\xi$ is as in \eqref{me1a}.
\end{lemma}

\bpf{Proof}\
We have for a.e. $(x,t) \in \Omega_\infty$ that
$$
|h_p(x,t)| \le H_0\left(|v(x,t)|^{p-1}+ 2|v(x,t)|^{p} + w_p(x,t)\right),
$$
where
\begin{eqnarray*}
w_p(x,t) &=& \frac{1}{p}\left(\int_0^t \xi(t-\tau)|v(x,\tau)| \dd\tau\right)^p\\
&\le& \frac{1}{p}\left(
\left(\int_0^t \xi(t-\tau) \dd\tau\right)^{1/p'}
\left(\int_0^t \xi(t-\tau)|v(x,\tau)|^p \dd\tau\right)^{1/p}\right)^p\\
&\le& \frac{1}{p} \int_0^t \xi(t-\tau)|v(x,\tau)|^p \dd\tau\,.
\end{eqnarray*}
Here, we have used H\"older's inequality with conjugate exponents $p,p'$
and Remark~\ref{remxi}. The assertion now follows
from Minkowski's inequality
$$
\left|\int_0^t \xi(t-\tau)|v(\cdot,\tau)|^p \dd\tau\right|_s \le
\int_0^t \xi(t-\tau)|v(\tau)|_{ps}^p \dd\tau\,.
$$
\epf

\begin{lemma}\label{ml3}
Let the hypotheses of Proposition \ref{moser} hold, and let $|\cdot|_{\infty,p}$
denote the norm in $L^\infty(0,\infty; L^p(\Omega))$. Let
$v \in W^{1,2}_{{\rm loc}}(0,\infty; (W^{1,2})'(\Omega))\cap
L^2_{{\rm loc}}(0,\infty;W^{1,2}(\Omega)) \cap L^\infty_{\rm loc}(\Omega \times (0,\infty))$
be a solution to (\ref{moser1})--(\ref{moser2}) such that $|v|_{\infty,p} < \infty$
for some $p\ge 1$. Then there exists a constant $\bar C>0$ independent of $v$ and $p$ such that
\begin{equation}\label{me8}
|v|_{\infty,2p} \le (\bar C p^{1+b})^{1/2p} \max\left\{1, V, V_\Gamma/C_f, |v|_{\infty,p}\right\}\,,
\end{equation}
where
\begin{equation}\label{me5}
b = \frac{N(q+1)}{2q - N}\,.
\end{equation}
\end{lemma}

\bpf{Proof}\
We test (\ref{moser1}) by $v |v|^{2p-2}$ to obtain, using Lemma \ref{ml2} that
\begin{eqnarray}\label{me2}
\frac{1}{2p} \frac{\dd}{\dd t} \io |v(x,t)|^{2p}\dd x
&+&\frac{2p-1}{p^2}\io |\nabla |v|^p|^2(x,t) \dd x +  \io |v(x,t)|^{2p}\dd x\\[2mm]\nonumber
&& + \, \ipo h(x)(C_f |v|^{2p} - v_\Gamma v |v|^{2p-2}\dd s(x)
\\[2mm]\nonumber
&\le& r^*H_0 \left(\frac{1}{2p} + 3 |v(t)|_{2pq'}^{2p} + \frac{1}{2p}
 \int_0^t \xi(t-\tau)|v(\tau)|_{2pq'}^{2p} \dd\tau\right).
\end{eqnarray}
We estimate the boundary integral using Young's inequality
$$
v_\Gamma v |v|^{2p-2} \le \frac{2p-1}{2p} C_f |v|^{2p} + \frac{1}{2p} C_f^{1-2p}V_\Gamma^{2p}\,.
$$
Set $v_p(x,t) = |v(x,t)|^p$. Then
\begin{eqnarray}\label{me3}
&&\hspace{-15mm}\frac{\dd}{\dd t} |v_p(t)|_2^2 + 2 |\nabla v_p(t)|_2^2 + 2p |v_p(t)|_2^2\\[2mm] \nonumber
&\le& C_f^{1-2p}|h|_{L^1(\partial\Omega)} V_\Gamma^{2p}
+ r^*H_0 \left(1 + 6p |v_p(t)|_{2q'}^{2} + 
 \int_0^t \xi(t-\tau)|v_p(\tau)|_{2q'}^{2} \dd\tau\right)\ \mae
\end{eqnarray}
By the Gagliardo-Nirenberg inequality, \cite{nier}, there exists a constant $G$
such that for all $\delta > 0$ and $t>0$ we have
\begin{equation}\label{me4}
|v_p(t)|_{2q'}^{2}\ \le\ G\left(\delta |\nabla v_p(t)|_2^2 + \delta^{-b}|v_p(t)|_1^2\right)
\ \le\ G\left(\delta |\nabla v_p(t)|_2^2 + \delta^{-b}|v|_{\infty,p}^{2p}\right),
\end{equation}
with $b$ given by (\ref{me5}). We now choose $\delta$ such that
$$
6pr^*H_0G\delta = 1\,,
$$
and obtain
\begin{eqnarray}\label{me6}
\frac{\dd}{\dd t} |v_p(t)|_2^2 + |\nabla v_p(t)|_2^2  + 2p |v_p(t)|_2^2 &\le&
C_7 \left(1 + (V_\Gamma/C_f)^{2p}+ p^{1+b} |v|_{\infty,p}^{2p}\right)\\[2mm]\nonumber
 &&+\, \frac{1}{6p} \int_0^t \xi(t-\tau)|\nabla v_p(\tau)|_{2}^{2} \dd\tau\quad \mae\,,
\end{eqnarray}
with a constant $C_7$ depending only on $\xi(0)$, $C_f$,
$|h|_{L^1(\partial\Omega)}$, $r^*$, $H_0$, and $G$. 
We now use Lemma \ref{ml1} with $x(t) = |v_p(t)|_2^2$, $y(t) = |\nabla v_p(t)|_2^2$,
$a = 2p$, $\delta = 1/(6p)$, $\mu \ge \min\{2, (5/6)\xi(0)\}$, and
$C = C_7 (1 + (V_\Gamma/C_f)^{2p}+ p^{1+b} |v|_{\infty,p}^{2p})$, which yields that
\begin{equation}\label{me7}
|v_p(t)|_2^2\ \le\ C_8\left(1 + V^{2p} + (V_\Gamma/C_f)^{2p} + p^{1+b}
|v|_{\infty,p}^{2p}\right)
\end{equation}
with a constant $C_8$ independent of $p$ and $t$, and (\ref{me8}) immediately follows.
\epf

We are now ready to finish the proof of Proposition \ref{moser}.

\bpf{Proof of Proposition \ref{moser}}\
For $k = 0, 1, 2, \dots$ set $y_k = \max\left\{1, V, V_\Gamma/C_f, |v|_{\infty,{2^k}}\right\}$.
We have $y_0 \le \max\left\{1, V, V_\Gamma/C_f, E_0\right\}$, and, as a consequence of
Lemma \ref{ml3},
$$
y_{k+1} \le (\bar C 2^{k(1+b)})^{2^{-(k+1)}} y_k\,.
$$
This yields
$$
\log y_{k+1} \le \log y_k + 2^{-(k+1)}(\log \bar C + k(1+b) \log 2)\,.
$$
Hence,
\begin{equation}\label{me9}
\log y_{n} \le \log y_0 + \sum_{k=1}^{n}\left(2^{-k}(\log \bar C + (k-1)(1+b) \log 2)\right)\,.
\end{equation}
The sum on the right hand side of (\ref{me9}) is convergent, and we easily complete the proof.
\epf

We now finish the proof of Theorem \ref{main} by showing that $T_0$ introduced at the beginning
of this subsection is $+\infty$ if $R$ is sufficiently large.
In (\ref{nequ2}), set again $U_{\Omega}(t) = \io U(x',t)\dd x'$. Then
$$
\dot U_{\Omega}(t)+(1+K_{\Gamma} |\Omega|)U_{\Omega}(t) = \io (\theta-\chi)(x',t)\dd x'
- |\Omega| p_0 \quad \mae
$$
By (\ref{esti1}), the right hand side of this ODE is uniformly bounded independently of $R$, hence
$|U_{\Omega}(t)| \le C_9$ in $(0,T_0)$. Using (\ref{nequ2}) once again, we obtain that
\begin{eqnarray}\label{me10}
|U(x,t)| &\le& C_{10}\left(1 + \int_0^t \expe^{\tau - t} \theta(x,\tau) \dd\tau\right)
\quad\mae\,,\\ \label{me11}
|U_t(x,t)| &\le& C_{11}\left(1 + \theta(x,t)
+ \int_0^t \expe^{\tau - t} \theta(x,\tau) \dd\tau\right)\quad\mae\,,
\end{eqnarray}
hence also (cf.~\eqref{nequ3})
\begin{equation}\label{me12}
|\chi_t(x,t)|\ \le\ C_{12}\left(1 + \theta(x,t)
+ \int_0^t\expe^{\tau - t}\theta(x,\tau)\dd\tau\right)\quad\mae
\end{equation}
As in (\ref{rhs}), we rewrite the right hand side of Eq.~(\ref{nequ1}) as
$$
- (\chi + U + p_0 + K_{\Gamma}U_{\Omega})U_t - (U+\chi + 1)\chi_t\,.
$$
By (\ref{esti1}), the function $U$ is in $L^\infty(0,\infty; L^2(\Omega))$ and the bound
does not depend on~$R$. Eq.~(\ref{nequ1}), with $\theta$ added to both the left and
the right hand side, thus satisfies the hypotheses of Proposition \ref{moser} for $N=3$ and $q=2$.
This enables us to conclude that $\theta(x,t)$ is uniformly bounded from above by a constant,
independently of $R$, so that $\theta$ never reaches the value $R$ if $R$ is sufficiently large,
which we wanted to prove. By (\ref{me10})--(\ref{me12}), also
$U$, $U_t$, and $\chi_t$ are uniformly bounded by a constant.

We proceed similarly to prove a uniform positive lower bound for $\theta$.
Set $R_0 := \sup\theta$, and in Eq.~(\ref{trequ1}) with $R>R_0$ put
$w = -\tilde w/\theta$, $\tilde w\in W^{1,2}(\Omega)$. For a new (nonnegative) variable $v(x,t) := \log R_0 - \log\theta(x,t)$
we obtain the equation
\begin{eqnarray}\label{me13}
&&\hspace{-15mm}\io v_t \tilde w(x)\dd x + \io \nabla v \cdot \nabla \tilde w(x) \dd x
+ \ipo h(x)\left(\frac{\theta_\Gamma}{\theta}-1\right) \tilde w(x) \dd s(x) \\ \nonumber
&=& \io\left(- \frac{U_t^2 +\chi_t^2}{\theta} -  \frac{|\nabla\theta|^2}{\theta^2}
+ U_t+ 2 \chi_t\right)  \tilde w(x)\dd x\,.
\end{eqnarray}
We now set
$$
\mathcal{H}[v] = \sign(v) \left(-\frac{U_t^2 +\chi_t^2}{\theta}-\frac{|\nabla\theta|^2}{\theta^2}
+ U_t+ 2 \chi_t\right)
$$
and check that the hypotheses of Proposition \ref{moser} are satisfied with
$f(v) = (\theta_\Gamma/R_0)(\expe^v - 1)$, $v_\Gamma = (R_0 -\theta_\Gamma)/R_0$,
$r \equiv 1$, and $v\mathcal{H}[v] \le 2C |v|$, where $C$ is a common upper bound for
$U_t$ and $\chi_t$. Hence, $v$ is bounded above by some $v^*$, which entails
$\theta \ge R_0 \expe^{-v^*}$.
This concludes the proof of Theorem \ref{main}.

\section{Long time behavior}\label{long}

In order to emphasize the relation to Section \ref{equi}, we keep the original
physical constants as in (\ref{sys2})--(\ref{sys3}). We prove the following statement.

\begin{proposition}\label{asp1}
Let the hypotheses of Theorem \ref{main} hold, and let the constants $\tilde\beta, \omega$
introduced in (\ref{1d17a}) satisfy the condition $1 -\tilde\beta > \max\{0,\omega\}$. Then we have
\begin{eqnarray}\label{ase1}
\int_0^\infty \left(\io \left(\theta_t^2 + U_t^2 + \chi_t^2 + |\nabla \theta|^2\right) \dd x
+ \ipo h(x)(\theta - \theta_\Gamma)^2\dd s(x)\right)\dd t &<& \infty\,,\\[2mm]\label{ase2}
\lim_{t\to\infty} \left(\io \left( U_t^2 + \chi_t^2 + |\nabla \theta|^2\right)(x,t) \dd x
+ \ipo h(x)(\theta - \theta_\Gamma)^2(x,t)\dd s(x)\right) &=& 0\,.
\end{eqnarray}
Furthermore, the functions $U_{\Omega}(t) = \io U(x,t) \dd x$,
$X_{\Omega}(t) = \io (1-\chi(x,t)) \dd x$ converge to their equilibrium values as $t \to \infty$.
\end{proposition}

The function on the left hand side of (\ref{ase2}) is almost everywhere equal
to a function of bounded variation. The limit is to be understood in this sense.

We see in particular that the temperature converges strongly in $W^{1,2}(\Omega)$
to its equilibrium value $\theta_\Gamma$, and the total ice contents $X_\Omega$ as well
as the pressure converge to a constant as $t\to \infty$.
For temperatures $\theta_\Gamma$ above the freezing point or below the undercooling limit,
this means, in agreement with the discussion in Section \ref{equi}, that also both $\chi(x,t)$
and $U(x,t)$ converge strongly in $L^1(\Omega)$ (hence, strongly in every $L^p(\Omega)$
for $p<\infty$) to their respective equilibrium values.
For intermediate temperatures, only the limit total ice contents can be identified,
but we are not able to decide about the convergence of the individual trajectories
to some of the equilibria.

\bpf{Proof}\
By (\ref{esti1}), since $\theta$ is uniformly bounded from above, we have
\begin{equation}\label{ase2bis}
\int_0^\infty\io \left(U_t^2 + \chi_t^2 + |\nabla \theta|^2\right)(x,t) \dd x \dd t
+ \int_0^\infty\ipo h(x)(\theta - \theta_\Gamma)^2(x,t)\dd s(x) \dd t < \infty\,.
\end{equation}
We rewrite Eq.~(\ref{sys2}) in the form
\begin{equation}\label{ase3}
c\theta_t - \kappa\Delta\theta = \nu (U_t)^2 - \beta\theta U_t + \gamma\chi_t^2 -
\frac{L}{\theta_c}\theta\chi_t\,.
\end{equation}
Due to the uniform $L^\infty$ upper bounds for $\theta$ and $U_t$, we can test Eq.~(\ref{ase3})
by $\theta_t$, and obtain
\begin{eqnarray}\label{ase4}
&&\hspace{-15mm}\io \theta_t^2(x,t) \dd x + \frac{\dd}{\dd t}\left(\io |\nabla \theta|^2(x,t) \dd x
+ \ipo h(x)(\theta - \theta_\Gamma)^2\dd s(x)\right)\\[2mm] \nonumber
&\le&  C_{13}\io \left(U_t^2 + \chi_t^2\right)\dd x\quad\mae\,,
\end{eqnarray}
with $C_{13}$ independent of $t$, and \eqref{ase2bis} with \eqref{ase4}
together with \cite[Lemma~3.1]{sz} yield
\begin{equation}\label{ase5}
\lim_{t\to\infty} \io |\nabla \theta|^2(x,t) \dd x
+ \ipo h(x)(\theta - \theta_\Gamma)^2(x,t)\dd s(x) = 0\,.
\end{equation}
System (\ref{sys1})--(\ref{sys3}) can be again considered as a the gradient flow of the form
(\ref{eg1}), with $v$ and $\psi(v)$ analogous to (\ref{eg4})--(\ref{eg5}), more precisely
\begin{eqnarray}\label{aseg4}
v &=& \vect{\nu U}{\gamma\chi}\,,\\[2mm]\label{aseg5}
 \quad \psi(v) &=& \io \left( \frac{\lambda}2 (U -\alpha(1-\chi))^2
+ (L\chi + \beta\theta_c U)\left(1-\frac{\theta_\Gamma}{\theta_c}\right)
+ I(\chi)\right)\dd x\\ \nonumber
&&+\, \frac{K_{\Gamma}}{2}\left(\io U \dd x + \frac{p_0}{K_{\Gamma}}\right)^2 + C_\psi\,,\\[2mm]\label{aseg6}
f &=& \vect{\beta(\theta - \theta_\Gamma)}{(L/\theta_c) (\theta - \theta_\Gamma)}\,.
\end{eqnarray}
We have $f \in L^2(0,\infty;L^2(\Omega)\times L^2(\Omega))$ by (\ref{ase2bis}), and
$\dot f \in L^2(0,\infty;L^2(\Omega)\times L^2(\Omega))$ by (\ref{ase4}).
{}From Lemma \ref{lg1} we conclude that
\begin{equation}\label{ase6}
\lim_{t\to\infty} \io \left(U_t^2 + \chi_t^2\right)(x,t) \dd x = 0\,.
\end{equation}
Set $\theta_{\Omega}(t) = \io \theta(x,t)\dd x$. The equation for $U_{\Omega}$ now reads
\begin{equation}\label{ase7}
\nu \dot U_{\Omega} + (\lambda+ K_{\Gamma}|\Omega|)U_{\Omega} = \alpha\lambda X_{\Omega}
+ \beta\theta_{\Omega} - (p_0 + \beta\theta_c)|\Omega|\,,
\end{equation}
hence
\begin{equation}\label{ase8}
\lim_{t\to\infty} \big((\lambda+ K_{\Gamma}|\Omega|)U_{\Omega}(t) - \alpha\lambda X_{\Omega}(t)\big)
= (\beta(\theta_\Gamma-\theta_c) - p_0)|\Omega|\,.
\end{equation}
{}From (\ref{sys1}) and (\ref{sys3}) we obtain for a.e. $(x,t) \in \Omega_\infty$ that
\begin{equation}\label{ase10}
\lambda(U - \alpha (1-\chi)) = -\nu U_t + \beta(\theta-\theta_c) - p_0 - K_{\Gamma} U_{\Omega}\,,
\end{equation}
\begin{equation}\label{ase11}
-\gamma \chi_t \in \alpha(-\nu U_t + \beta(\theta_\Gamma-\theta_c)
- p_0 - K_{\Gamma} U_{\Omega}) + L\left(1- \frac{\theta_\Gamma}{\theta_c}\right)
+ \left(\alpha\beta - \frac{L}{\theta_c}\right)(\theta - \theta_\Gamma) + \partial I(\chi)\,.
\end{equation}
We define an auxiliary function
\begin{eqnarray}\label{ase9}
A(x,t) &:=& -\gamma\chi_t(x,t) +\alpha \nu U_t(x,t) - \left(\alpha\beta - \frac{L}{\theta_c}\right)
(\theta(x,t)-\theta_\Gamma)\\[2mm]
&& +\, \alpha K_{\Gamma} U_{\Omega}(t)
- \frac{\alpha^2\lambda K_{\Gamma} X_{\Omega}(t)
+ \alpha(\beta (\theta_\Gamma-\theta_c) - p_0)K_{\Gamma}|\Omega|}
{\lambda+ K_{\Gamma}|\Omega|}\,.
\end{eqnarray}
With the notation (\ref{1d17}), (\ref{1d17a}) we rewrite (\ref{ase11}) in the form
\begin{equation}\label{ase12}
\frac{1}{L}A(x,t) + \frac{d}{|\Omega|} X_{\Omega}(t) + (1-\tilde\beta) 
\left(\frac{\theta_\Gamma}{\theta_c} - 1\right)
+ \omega \in \partial I(\chi(x,t))\quad \mbox{a.e.}\,,
\end{equation}
as an evolution counterpart of the equilibrium condition (\ref{1d15a}).
The above computations show that $\lim_{t\to\infty}|A(t)|_2 = 0$. We now prove the following
implications:
\begin{itemize}
\item[{\rm (i)}] If $\theta_\Gamma \ge \theta_c(1-\omega/(1-\tilde\beta))$ then
$\io|1 - \chi(x,t)|\dd x \to 0$ as $t\to \infty$;
\item[{\rm (ii)}] If $\theta_\Gamma
\le \theta_c(1-(\omega + d)/(1-\tilde\beta))$ then
$\io \chi(x,t) \dd x \to 0$ as $t\to \infty$;
\item[{\rm (iii)}] If $\theta_c(1-(\omega + d)/(1-\tilde\beta))<\theta_\Gamma
< \theta_c(1-\omega/(1-\tilde\beta))$ then\\
$X_{\Omega}(t) \to (|\Omega|/d)((1-\tilde\beta)(1 - (\theta_\Gamma/\theta_c))-\omega)$
 as $t\to \infty$.
\end{itemize}
The corresponding convergence of $U$ then follows from (\ref{ase8})--(\ref{ase10}).

To prove the above statements (i)--(iii), set
$$
\chi^* := \frac{1}{d} \left(
(1-\tilde\beta)\left(1-\frac{\theta_\Gamma}{\theta_c}\right) -\omega \right)\,,
\quad A^*(x,t) := \frac{1}{Ld} A(x,t)\,.
$$
Eq.~(\ref{ase12}) reads
\begin{equation}\label{ase13}
A^*(x,t) + \frac{1}{|\Omega|} X_{\Omega}(t) - \chi^* \in \partial I(\chi(x,t))\quad \mbox{a.e.}\,,
\end{equation}
that is,
\begin{equation}\label{ase14}
\left(A^*(x,t) + \frac{1}{|\Omega|} X_{\Omega}(t) - \chi^*\right)(\tilde \chi - \chi(x,t)) \le 0
\quad \mbox{a.e.} \ \ \forall \tilde\chi \in [0,1]\,.
\end{equation}
Integrating over $\Omega$, we obtain for every $\tilde\chi \in [0,1]$ and a.e. $t>0$ that
\begin{equation}\label{ase15}
\left(\frac{1}{|\Omega|} X_{\Omega}(t) - \chi^*\right)
\left(\frac{1}{|\Omega|} X_{\Omega}(t) - (1-\tilde\chi)\right)\le -\frac{1}{|\Omega|}
\io A^*(x,t)(\tilde \chi - \chi(x,t)) \dd x\,.
\end{equation}
The right hand side of (\ref{ase15}) tends to $0$ as $t$ tends to $\infty$.
Hence,
\begin{equation}\label{ase16}
\limsup_{t\to\infty} \left(\frac{1}{|\Omega|} X_{\Omega}(t) - \chi^*\right)
\left(\frac{1}{|\Omega|} X_{\Omega}(t) - (1-\tilde\chi)\right) \le 0
\quad \forall \tilde\chi \in [0,1]\,.
\end{equation}

\begin{itemize}
\item[{(i)}] We have $\chi^* \le 0$. The assertion follows if we put
$\tilde \chi = 1$ in (\ref{ase16}).
\item[{(ii)}] We have $\chi^* \ge 1$. The argument of (i) applies if we put
$\tilde \chi = 0$ in (\ref{ase16}).
\item[{(iii)}] Here, we have $0 < \chi^* < 1$, and it suffices to put
$\tilde\chi = 1 - \chi^*$.
\end{itemize}
\epf

Note that in all cases the difference $U - \alpha(1-\chi)$ converges in $L^2(\Omega)$
to its equilibrium value as $t\to \infty$. The problem if $\chi(x,t)$ and $U(x,t)$
separately converge in the case (iii) is still open.

\noindent {\bf Acknowledgments}.
The authors wish to appreciate helpful suggestions and comments by Ch\'erif Amrouche
and Ji\v r\'{\i} Neustupa.

\end{document}